\documentclass[journal=jacsat,manuscript=article]{achemso}
\usepackage[utf8]{inputenc}
\usepackage[compact]{titlesec}
\title{A semi-continuous model for transmission of SARS-CoV-2 and other respiratory viruses in enclosed spaces via multiple pathways to assess risk of infection and mitigation strategies}

\author{Panagiotis Demis}
\affiliation[University of Surrey]
{Department of Chemical and Process Engineering, University of Surrey, Guildford, UK}

\author{Ishanki De Mel}
\affiliation[University of Surrey]
{Department of Chemical and Process Engineering, University of Surrey, Guildford, UK}

\author{Hayley Wragg}
\affiliation[University of Bath]
{Department of Mathematical Sciences, University of Bath, Bath, UK}

\author{Michael Short}
\affiliation[University of Surrey]
{Department of Chemical and Process Engineering, University of Surrey, Guildford, UK}

\author{Oleksiy V. Klymenko}
\affiliation[University of Surrey]
{Department of Chemical and Process Engineering, University of Surrey, Guildford, UK}
\email{o.klymenko@surrey.ac.uk}

\date{\today}
\usepackage[titletoc,title]{appendix}
\usepackage{Packages}
\usepackage{Notation}
\usepackage{Figures}

\begin{document}
\maketitle

\section{Abstract}
The Covid-19 pandemic has taken millions of lives, demonstrating the tragedy and disruption of respiratory diseases, and how difficult they can be to manage. However, there is still significant debate in the scientific community as to which transmission pathways are most significant and how settings and behaviour affect risk of infection, which all have implications for which mitigation strategies are most effective. This study presents a general model to estimate the rate of viral transfer between individuals, objects, and the air. The risk of infection to individuals in a setting is then computed considering the behaviour and interactions of individuals between themselves and the environment in the setting, survival times of the virus on different surface types and in the air, and mitigating interventions (ventilation,  hand disinfection, surface cleaning, etc.). The model includes discrete events such as touch events, individuals entering/leaving the setting, and cleaning events. 

We demonstrate the model capabilities on three case studies to quantify and understand the relative risk associated with the different transmission pathways and the effectiveness of mitigation strategies in different settings. The results show the importance of considering all transmission pathways and their interactions, with each scenario displaying different dominant pathways depending on the setting and behaviours of individuals therein.

The flexible model, which is freely available, can be used to quickly simulate the spread of any respiratory virus via the modelled transmission pathways and the efficacy of potential mitigation strategies in any enclosed setting by making reasonable assumptions regarding the behaviour of its occupants. It is hoped that the model can be used to inform sensible decision-making regarding viral infection mitigations that are targeted to specific settings and pathogens.

\section{Introduction}
After the sudden emergence and spread of the novel coronavirus, SARS-CoV-2, the World Health Organisation (WHO) declared a pandemic on 11 March 2020. The COVID-19 pandemic has demonstrated the global impact, tragedy, and disruption caused by viral outbreaks, and the scientific community has warned that it is unlikely to be the last such outbreak. Despite the incredible efforts of the world's governments and scientific community, there is still significant debate about which viral transmission pathways contribute most significantly to the spread of SARS-CoV-2 between individuals, and correctly characterising this has important repercussions for decision-making regarding mitigation strategies and what economic activities should be allowed to continue during pandemics \citep{RAMPvent21}.

SARS-CoV-2 and other respiratory viruses are thought to be spread through 3 major transmission mechanisms: direct and indirect respiratory transfer via respiratory droplets and aerosols, and through the direct contact pathway \citep{mittal_ni_seo_2020}. Direct respiratory transfer takes place when an infected individual secretes large droplets ($>$ 5 $\mu$m) via vocalising, coughing, sneezing etc. that are directly deposited onto susceptible individuals' mucous membranes via close proximity. Early in the COVID-19 pandemic, this was thought to be the main transmission pathway and this is mitigated quite effectively through the use of facemasks \citep{morawska-dropletsizes, Grover2020mask}.

Since, the aerosol/microdroplet pathway has also been identified as a potential respiratory transmission pathway \citep{who-airborne}. In this pathway, small airborne particles secreted by infected individuals remain suspended in the air and can be inhaled by susceptible individuals. This has been observed for other viruses such as influenza A \citep{reviewoffluviaaerosol} and has also been conjectured to have been responsible for SARS-CoV-2 transmission events, such as in a Chinese restaurant  \citep{Li2020.04.16.20067728}, a cruise ship \citep{cruiseship21}, and the Skagit Valley Chorale superspreading event \citep{skagit}. Mitigation of this pathway is challenging, but increasing ventilation rates with more outside air, discouraging high concentrations of people indoors and limiting the time spent indoors are most effective \citep{morawska20-airbornemin, RAMPvent21}. It has also been shown in \citet{leungrespiratory} that face coverings can reduce aerosol emissions. However, other studies have found that this effect can be limited due to small aerosol particles either penetrating the mask or leaking through gaps around the cheeks and nose, which results in apparent filtration efficiencies of face coverings from cloth masks to N95-grade respirators in the range of 10-60\% \cite{Bundgaard2020, Shah2021}.

Lastly, the contact route is a mechanism in which the surface of an object is contaminated with viral particles, either via deposition from the air (through particle settling or droplets) or from an individuals' contaminated hand \cite{XiaoRoleofFomites}. Such a contaminated object is called a fomite. Susceptible individuals can then transfer the virus from the fomite to their mucous membranes through touching the object and then their face \cite{Kraayfomite2018}. This mechanism has been anecdotally recognised by several hospitals \citep{SAHospitalReport} and the WHO and UK governmental advice to frequently wash hands and disinfect surfaces shows that this mechanism is believed to be important. \Cref{fig::pathway} represents the transmission pathways for viral infections.

\Pathway

Most epidemiological studies related to viral disease transmission have focused on macroscopic, large-scale transmission and often use modified classical susceptible-infected-recovered (SIR) modelling frameworks \cite{miao20,Wijaya2020.04.25.20079178}. There have been some efforts to use these models to determine which transmission routes are likely to be dominant\citep{Zhang14857,Ferrettieabb6936}, however they are not able to consider environment-specific risks. Many researchers have also performed detailed transmission studies and simulations to understand and quantify the specific aerosol \cite{morawska-airbornetrans, burridge2020airborne, YouAirline}, droplet \citep{NanZhangBehaviour, Zhang2020-dominantroutes,mittal_ni_seo_2020}, and contact routes \citep{Kraay2020.08.10.20171629, king2020fingertips}, however few modelling efforts exist for quantifying overall risks of infection based on the combination of transmission routes.

Fomite transmission routes have seen significant study via models and simulations. \citet{Zhao2012} developed an Environmental Infection Transmission System model to quantify risks of infection from droplet-contaminated and hand-contaminated routes, concluding that public, large-surface area droplet-contaminated surfaces have the highest transmission potential. \citet{Beamer2015} used micro-activity empirical data to validate their model for fomite pathogen transmission to test different workplace strategies for reducing infection risks from rhinovirus and rotavirus. The effects of increasing hand hygiene to reduce potential infection by SARS-CoV-2 was investigated by \citet{Pham2020handhygiene} and their results show that event-based hand washing (such as after touching an object) is more effective than frequency washing (every 30 minutes, for example). \citet{Kraayfomite2018} developed an ordinary differential equation (ODE) model for fomite-mediated virus transmission, with results highlighting that fomites play an important role in virus transmission. The study also simulated cleaning events on hands and surfaces, demonstrating the efficacy of different cleaning strategies. 

To understand influenza transmission routes, \citet{xiao_tang_hui_lei_yu_li_2018} used a multi-agent modelling framework to understand a nosocomial outbreak in a Hong Kong hospital. Their detailed spatio-temporal model identified that long-range airborne (94 \%) and fomite routes (6 \%) were both likely to have played a role in the outbreak. \citet{xiao2017norovirus} used a similar multi-agent approach to understand the transmission of a norovirus outbreak in a UK hotel restaurant in 1998 and found that, out of the multiple pathways examined, fomite transmission played the largest role in the outbreak.

\citet{Lei2018} developed a multi-route disease transmission model to study in-flight outbreaks of norovirus, SARS-CoV and influenza A H1N1. They used the detailed seat positions of those infected to consider all 3 routes, modelling each of the 3 transmission routes separately. Markov Chains were built to model fomite transmission routes in a surface contamination network via transfer efficiencies and transition matrices. Their study concluded that for H1N1 transmission, the close contact and aerosol routes were more important, but for SARS-CoV, fomites are slightly more important than the other transmission routes, although all 3 play a key role. Finally for norovirus, fomites are the largest contributor to infection risk. Using a similar Markov Chain modelling approaches, \citet{cruiseship21} showed that it was likely that multiple transmission routes played a role in the SARS-CoV2 outbreak on the Diamond Princess cruise ship. \citep{relpathwaysjones} used a similar approach to attempt to quantify relative contributions of different pathways in healthcare personnel in patient care, determining that all pathways played a part, however, respiratory pathways dominated during the short patient interactions. These studies highlight the importance of developing models to understand the various transmission pathways, quantify the importance of each route, and the environment-specific risks.

\citet{ZhangLi2018transmission} developed a model considering the three different routes of transmission for the Influenza A virus in a student office (see \Cref{fig::pathway}). They used observations from camera recordings to track interactions between students and model the risk of infections to individuals. \citet{ZhangLi2018transmission} present two equations to determine the total quantity of virus on a hand $\VirHand $ and on a surface $\VirSurf $ at time $t$:

\begin{equation} 
\label{eqn::DivVirHand}
 \dv{\VirHand}{t}=\TraSurfHand{} \AreaCon{} \frac{\VirSurf}{ \Area[\surf] } 
                    - \TraHandSurf{} \AreaCon{}\frac{ \VirHand}{ \Area[\hand]} 
\end{equation}
\begin{equation} 
\label{eqn::DivVirSurf} 
\dv{\VirSurf}{t}=\TraHandSurf{} \AreaCon{}\frac{ \VirHand}{ \Area[\hand]}
                    -\TraSurfHand{} \AreaCon{} \frac{\VirSurf}{ \Area[\surf]} 
\end{equation}

where $\TraSurfHand{}$ is the transfer rate from surface to hand; $\TraHandSurf{}$ represents the transfer rate from hand to surface; $\AreaCon{}$ is the area of the surface that made contact; $\Area[\hand]$ and $\Area[\surf]$ represent areas of the hand and surface, respectively. Through a detailed surface touch map, found via video surveillance, they found that both public surfaces and private surfaces play a part in spreading viruses via fomites \citep{zhang2018surface}. They modeled close contacts by observing the number of contacts and duration, and making assumptions regarding the transmission of small and large droplets between individuals. There results suggest that, for influenza A, 54.3 \% of infection risk is attributed to the aerosol pathway, 44.5\% to the droplet pathway and 4.2\% to the fomite pathway, demonstrating the importance of considering multiple pathways in viral transmission modelling. In subsequent work, \citet{zhang20closecont} the authors obtained more accurate parameter values for close contacts. Their approach of combining data with different approximate models is promising and effective, however, up to now the different routes were modelled separately and the model has yet to be generalised to scenarios beyond those for which data is directly available. 

A recent review from \citet{LeungNature21} highlights the lack of understanding and intense debate surrounding the importance of the different transmission mechanisms for various respiratory diseases, stressing the need for inter-disciplinary cooperation to quantify and understand these mechanisms holistically. Relatively few models attempt to quantify the effects of multiple transmission routes simultaneously, and fewer still can be used to simulate and understand the effects of different mitigation strategies. While many models with detailed spatio-temporal distributions exist using computational fluid dynamics, multi-agent models, etc., these are limited in their applicability to general scenarios, are computationally costly to run, and require detailed mathematical and coding knowledge to test different scenarios. On the other hand, models that only consider population-level transmission, can be difficult to assess for local decision-making and contain a number of highly inaccurate parameters and many simplifications.

This study presents a novel, general and flexible model for risk of viral infection in enclosed spaces. The model, which has been applied to the SARS-CoV-2 virus in our case studies, can be used to study any viral transmission by modifying the appropriate virus-specific parameters and used for modelling any enclosed environment (such as public transport, offices, and classrooms) by modifying the environmental parameters. It can be used to approximate the risk of infection of individuals, based on the natural decay of viruses on different surface types and in the air, the surface types, transfer efficiencies of the virus from the surfaces, ventilation rates, and human behaviour in the settings (interactions with objects, face touching frequencies, etc.). The deterministic model is formulated as a set of first-order ODEs, that also has discrete inputs that are used to model events such as individuals entering or leaving the enclosed environment and cleaning events. The model can be used to quantify and understand the relative risk associated with the different transmission pathways and the effectiveness of different mitigation strategies in different environments.

\section{Mathematical model \label{sec::model}}

In the following, we will be considering a generic enclosed setting such as an office, a classroom, a retail outlet, a gym or public transport (e.g., train carriage, bus, etc.) which individuals can enter at different times and remain therein for different periods of time. The way individuals interact with each other and objects in the setting is parameterised in the model since it is dependent on setting type and associated behavioural patterns. For example, in an office environment each individual interacts primarily with few objects on and around their desk (‘private surfaces’ following the terminology of \citet{ZhangLi2018transmission}) and only occasionally touches shared objects (or ‘public surfaces’ \cite{ZhangLi2018transmission}) such as door handles, light switches, water fountains, etc. On the other hand, people visiting retail outlets or gyms or using public transport interact mainly with shared (public) objects, and the associated touching behaviour is markedly different from that in a typical office.

The model tracks the time-dependent spread of virus in the setting from one or more ‘sources’ (infected individuals) to mucous membranes and respiratory tracts of susceptible individuals through the transmission routes shown in \Cref{fig::pathway} while accounting for the effects of natural deactivation of the virus in air and on fomite surfaces \cite{van2020aerosol} as well as a range of mitigation interventions. The latter can be either preventative (e.g. the wearing of face coverings) or mitigatory (physical removal of the virus through, e.g., ventilation, air filtration, cleaning and handwashing, or viral deactivation via disinfection, sanitising, UV irradiation, etc.).

Consider an enclosed space (setting) characterised by volume $\volair$ that is visited by $\numhands$ individuals over a period of $\tend$ hours. Individuals can arrive at different times $\tin[j]$ and remain in the setting for the duration $\tdur[j]$. The presence of an individual in a setting can be described by the following indicator function:
\begin{equation}
    \label{eqn::IndicHeaviside}
    \Indic[j](t)=\Theta(t-\tin[j]) \times \Theta(\tin[j]+\tdur[j]-t), \ \ j = 1,\ldots,\numhands
\end{equation}
where $\Theta(\cdot)$ is the Heaviside step function.\
This formalism allows modelling people getting on and off public transport, entering and leaving shops, offices and other enclosed settings, and can be easily extended to multiple visits by an individual to the same setting.

The status of an individual is represented by a binary parameter $\infec[j], \ j = 1,\ldots,\numhands$ which takes the value 1 if they are infected or 0 if they are susceptible:
\begin{equation}
    \infec[j]=\left\lbrace \begin{array}{cc}
     1    & \text{if } j \text{ is infected} \\
     0    & \text{otherwise}
    \end{array} \right.
\end{equation}

We assume that at the time of entry susceptible individuals are free of the virus, i.e., their viral loads on the hands, $\VirHand[j](t)$, and mucous membranes, $\VirMems[j](t)$ are zero at $t=0$, while infected individuals carry a significant amount of the virus on their mucous membranes that is assumed to remain constant throughout the period of interest, and the viral load on their hands at $t=0$ is assumed to be in equilibrium with their mucous membranes (see \cref{sec::FullModel}).

The setting is assumed to contain $\numsurf$ objects that can be touched or handled by people present in the setting leading to their contamination, the level of which is characterised by viral loads $\VirSurf[i](t), \ i = 1,\ldots, \numsurf$. Surface contamination can also occur through the deposition of viral particles contained in droplets/aerosol expelled by infected individuals. Contamination of air in the setting by virus-laden aerosol is characterised by the overall viral load in the air, $\VirAir(t)$. A susceptible individual can also pick up virus from one surface with their hands and transfer it to another, which leads to the spread of the virus through the environment.

The purpose of the model is to quantify the viral load, $\VirMems[j](t)$, accumulated by a susceptible individual $j$ via all three transmission pathways while in the setting as well as their associated risk of developing COVID-19:
\begin{equation}
\label{eqn::InfectionRisk}
\infecP[j]=1-\exp\left(-\frac{\VirMems[j](t)}{\doserep[j]}\right)
\end{equation}
which depends on the dose response parameter $\doserep[j]$.

We note that the latter measure may not be reliable due to the high uncertainty in the value of the dose response parameter, $\doserep[j]$, measured in the units of $\left(\text{virus quantity}^{-1}\right)$. Therefore both the risk of infection and the cumulative viral load, $\VirMems[j]$, are considered as important outputs of the model presented.

The above-mentioned transmission mechanisms along with the underlying modelling assumptions are described mathematically below. In the following subsections we are developing a continuous model formulation relying on average frequencies of discrete events such as touching fomite surfaces, breathing, coughing, etc., unless otherwise stated.

\subsection{Aerosol pathway \label{sec::aeropath}}
During their stay in the setting, individuals are assumed to expel droplets of mucus of various sizes via their respiratory activities (including breathing, coughing and sneezing) and vocalisation. The corresponding rate of continuous viral shedding through these activities by individual $j$ is denoted in the following as $\airrate[j]$. While this rate should clearly be dependent on the viral load on an individual's mucous membranes \cite{Kawasuji2020virload} the form of this dependence is unknown. There are, however, experimental measurements of the rate of virus shedding by infected individuals through this pathway \cite{Ma2020shedding}. Therefore, $\airrate[j]$ is treated as a parameter in our model, the values of which are presented in \Cref{tab:table_virus} of Appendix \ref{sec::AppendixA}.

Following the widely-adopted terminology, we will classify droplets expelled through respiratory activities and vocalisation into large and small \cite{ZhangLi2018transmission, Lei2018}. This subdivision is based on the propensity of small droplets to remain aerosolised for prolonged periods of time in that they are either small enough when expelled or rapidly lose their water content through evaporation before hitting the ground leaving aerosolised nuclei. By contrast, large droplets do not reduce in size sufficiently so as to become aerosolised and land on nearby surfaces, which may lead to their contamination.

In the model, the fraction of large droplets expelled by individual $j$ is denoted as $\fractionLD[j]$. These droplets are involved in the close contact transmission pathway which will be discussed in Section \ref{sec::LDpath}. 

Using this notation, the rate at which an individual $j$ expels small droplets, which contribute to the viral load suspended in the air and can subsequently be inhaled by other individuals, can be expressed as

\begin{equation}
\left(1-\fractionLD[j]\right)\airrate[j]
\end{equation}

Note that we do not consider individual respiratory events, so the shedding of droplets of mucus is treated as a continuous process for the purpose of modelling.

We assume that the air in the setting is perfectly mixed such that aerosol generated anywhere in the room is instantly dispersed equally throughout the room. This approximation is often used in modelling aerosol distribution in enclosed settings since the time scale of air homogenisation is often significantly shorter than other time scales, particularly the duration of time individuals spend in the setting.

\par
The deposition of small virus-laden droplets/aerosols onto fomite surfaces as well as on human skin has also been included in the model. The dynamics between airborne and fomite virus concentrations are discussed in \citet{ZhangLi2018transmission}, including loss of the virus from surfaces through resuspension. The rate of virus deposition onto surface $i$ is proportional to its surface area and the viral concentration in air so that

\begin{equation}
\SDrate[i]= \deposconst \Area[{\surf[i]}]\frac{\VirAir}{\volair}
\end{equation}
where $\deposconst$ is the rate constant for small droplet deposition.
\par
While we include in our model, for the sake of generality, the terms describing the rates of deposition of small droplets, $\SDrate[i]$, and resuspension of the virus from surface $i$ back into the air, $\resusprate[i]$, there is no sufficient evidence in the literature to quantify these phenomena. Thus, although the corresponding terms are included in the formulation, they are kept zero in all the simulations reported below.

\subsection{Close contact pathway \label{sec::LDpath}}
Large droplets expelled by an infected individual may travel only a short distance from the source before landing on nearby surfaces, other people's skin, or mucous membranes. The viral load deposited on an object through this pathway depends on the following:
\begin{itemize}
\item The rate at which large droplets are emitted by the infected individual, $\fractionLD[j]\airrate[j]$ (see Section \ref{sec::aeropath}).
\item The amount of time the infected individual $j$ spends in close proximity to the object or other individual $x$. This is represented in the model by the fraction of the time, $\fractionCloseContact[j]{x}$, that the pairs spend in close proximity.
\item The fraction of large droplets emitted by the infected individual $j$ that land on object $x$, $\fractionLDsurf[j]{x}$, which depends on the relative positions and the distance between the source of the droplets and the acceptor surface. \cite{zhang20closecont} The latter can also change with time, so this parameter should be interpreted as the average fraction of large droplets transferred while in close proximity.
\end{itemize}
Thus, the rate of viral deposition from source $j$ onto surface $x$ through the close contact pathway is given by
\begin{equation}
\LDrate[j]{x} = \fractionCloseContact[j]{x}\fractionLDsurf[j]{x}\fractionLD[j]\airrate[j]
\end{equation}
while both the individual $j$ and acceptor $x$ are present in the setting. As with the aerosol generation, the deposition of large droplets is treated as a continuous process, and the above rate can be considered as the average rate of large droplet deposition throughout the duration of close contact between the individual $j$ and acceptor surface $x$.

Viral particles deposited onto surfaces with large droplets can then be spread further by individuals through touch (fomite-mediated transmission), which is discussed in \ref{sec::fomitepath}.

\subsection{Fomite transmission pathway \label{sec::fomitepath}}

The fomite transmission pathway involves the transfer of viral particles between an individual’s hand and a fomite surface or their own mucous membranes (mouth, nose or eyes). Through this pathway, an infected individual can contaminate their hands by touching their mouth or nose before handling other objects while all individuals can spread the virus by touching multiple objects. Susceptible individuals can self-inoculate by transferring the virus from their contaminated hands to mucous membranes.

\par \Cref{fig::surfacepath} shows the possible fomite transmission pathways between surfaces, a hand and a mucous membrane as well as the deposition of virus-laden droplets and aerosol onto surfaces (see Sections \ref{sec::aeropath} and \ref{sec::LDpath}).

\SurfacePath

In the following, we will be treating the surface touch behaviour as a continuous process as opposed to individual touch events. We assume that the rate of transfer from a donor object \(x\in\left[\surf[1],...,\surf[\numsurf],\hand[1],...,\hand[\numhands], \mems[1],...,\mems[\numhands] \right]\) to an acceptor object $y\in[\surf[1],...,\surf[\numsurf],\hand[1],...,\hand[\numhands], \\   \mems[1],...,\mems[\numhands]]$ (i.e., these can be surfaces, hands, or mucous membranes) is proportional to the average concentration of the virus on $x$, $\Vir[x]/\Area[x]$, average frequency of contacts between $x$ and $y$, $\frequency[x]{y}$, and contact area, $\AreaCon[x]{y}$, specific to the way objects $x$ and $y$ come into contact. $\Vir[x]$ is the viral load on $x$ and $\Area[x]$ is its surface area. Using this notation, the rate $\transferrate[x]{y}$ of viral transfer from object $x$ to object $y$ can be described using the following formalism:
\begin{equation}
\label{eqn::RateXY}
\transferrate[x]{y}=\fractionrate[x]{y}\frequency[x]{y}\AreaCon[x]{y} \frac{\Vir[x]}{\Area[x]}=\transferconst[x]{y}\Vir[x]
\end{equation}
where \(\fractionrate[x]{y}\) is the fraction of the viral load on surface $x$ within the area of contact $\AreaCon[x]{y}$ that is transferred to $y$ upon one contact, and $\transferconst[x]{y}=\fractionrate[x]{y}\frequency[x]{y}\AreaCon[x]{y}/\Area[x]$ is the overall transfer rate constant from $x$ to $y$. Note that $\frequency[x]{y}=\frequency[y]{x}$ and $\AreaCon[x]{y}=\AreaCon[y]{x}$ for any $x$ and $y$.

\par In our case studies presented in Section \ref{sec::casestudies} we assume that there are no direct contacts between individuals (i.e., the hands and mucous membranes of one individual do not come into contact with either hands or mucous membranes of another), which corresponds to physical distancing being enforced. Therefore, $\fractionrate[H_i]{M_j}=\fractionrate[M_i]{H_j}=\fractionrate[M_i]{M_j}=0$ if $i\neq j$.
Direct contacts between fomites are also not considered, so $\fractionrate[S_i]{S_j}=0$ for any $i$ and $j$.

\subsection{Removal and inactivation of the virus \label{sec::removal}}

\subsubsection{Continuous formulation}
\par 
Natural inactivation of the virus is assumed to be exponential in air and on all surfaces, including skin, so that the absolute rates of inactivation can be described as
\begin{equation}
\label{eqn::Inactivation}
\inact[x]=\inactconst[x]\Vir[x]
\end{equation}
where $\inactconst[x]$ is the inactivation rate constant which can be defined through virus half-life, $\halflife[x]$, on surface $x$ or in air as $\inactconst[x]=\ln(2)/\halflife[x]$. The latter parameter is dependent on the environmental conditions (temperature and relative humidity) and, for fomites, on their material type (e.g., stainless steel, copper, plastic, paper, etc.) and structure of the surface (porous or non-porous). These dependencies can be readily accounted for in the model using an Arrhenius-type inactivation model such as the one proposed by \citet{Yap2020inactivation}.

\par 
In the continuous formulation, the rates of any washing/cleaning interventions are assumed to be described by continuous rates which, on average, correspond to the frequency of cleaning events.
\par
The absolute rate of decrease of viral load on the hands due to washing (or using hand sanitiser gel, wipes, etc.) can be described as:
\begin{equation}
\label{eqn::HandWashRate}
\TraCleanHand[j]{j}(t) =\cleaneffecHand[j]\freqcleanHand[j]\VirHand[j]
\end{equation}
where $\cleaneffecHand[j]$ is the ‘efficiency’ of virus removal during hand-washing, and $\freqcleanHand[j]$ is the frequency of hand washing. If the fraction of the virus removed in a hand-washing event, $\cleanfracHand[j]$, is given the two parameters are related through \footnote{\label{fn_washrate} The rate of viral load reduction due to hand washing expressed using the continuous formulation is
\begin{equation}
\dv{\VirHand[j](t)}{t} =-\TraCleanHand[j]{j}=-\cleaneffecHand[j]\freqcleanHand[j]\VirHand[j]
\end{equation}
which describes an exponential decay
\begin{equation}
\VirHand[j](t) =\VirHand[j](t=0) \, e^{-\cleaneffecHand[j]\freqcleanHand[j]t}
\end{equation}
The relative reduction of the viral load between hand washing events is then
\begin{equation}
\frac{\VirHand[j]\left(t=1/\freqcleanHand[j]\right)}{\VirHand[j](t=0)} = e^{-\cleaneffecHand[j]} =\cleanfracHand[j]^{-1}=10^{-\lrv[j]}
\end{equation}
}

\begin{equation}
\cleaneffecHand[j]=\ln(\cleanfracHand[j])
\end{equation}
and if $\cleanfracHand[j]$ is represented using the $\log10$ reduction value ($\lrv$), $\cleanfracHand[j]=10^{\lrv[j]}$,
\begin{equation}
    \cleaneffecHand[j]=\lrv[j]\ln(10)
\end{equation}

A similar expression can be used for the cleaning of fomite surfaces:
\begin{equation}
    \TraCleanSurf[i]=\cleaneffecSurf[i]\freqcleanSurf[i]\VirSurf[i]
\end{equation}
where $\cleaneffecSurf[i]$ is the virus deactivation efficiency on surface $i$, and $\freqcleanSurf[i]$ is the frequency of cleaning the surface. If the fraction of the virus removed in a single intervention, $\cleanfracSurf[i]$, is given, the two parameters are related through
\begin{equation}
    \cleaneffecSurf[i]=\ln(\cleanfracSurf[i])
\end{equation}
and if $\cleanfracSurf[i]$ is represented using $\log10$ reduction value ($\lrv$), $\cleanfracSurf[i]=10^{\lrv[i]}$,
\begin{equation}
\cleaneffecSurf[i]=\lrv[i]\ln(10)
\end{equation}

\subsubsection{Mixed continuous-discrete formulation}
Some or all of the events such as touching surfaces, washing hands and cleaning of surfaces can be modelled as discrete occurrences. This can be combined with continuous formulations for other transmission rates, however the form of some of the terms would change.

If the frequency of contacts between the hands of individuals and various surfaces is much higher than the rate of handwashing and cleaning, it is justified to describe only the latter two using a discrete formulation.
Hand washing (or using hand sanitiser gel, wipes, etc.) at discrete times $t^w_k$ is assumed to lead to an instantaneous drop in virus concentration on the hands by a factor $\cleanfracHand[j]$ (see the footnote on page \pageref{fn_washrate}) so that the absolute rate of viral removal is:
\begin{equation}
\label{eqn::handcleanrate}
\TraCleanHand[j]{}=\sum_{k=1}^{\numhwash[j]}\cleanfracHand[j]\VirHand[j]\delta(t^w_k)
\end{equation}
where $\numhwash[j]$ is the total number of times individual $j$ washes their hands, and $\delta$ is the Dirac delta-function.
\par
A similar expression can be used for the cleaning of fomite surfaces:
\begin{equation}
\label{eqn::fomitecleanrate}							
    \TraCleanSurf[i]{}=\sum_{k=1}^{\numclean[i]}\cleanfracSurf[i]\VirSurf[i]\delta(t^c_k)
\end{equation}
where $\numclean[i]$ is the number of times surface $i$ is cleaned, and $\cleanfracSurf[i]$ is the fraction of the virus deactivated by each of the interventions taking place at times $t^c_k$.

\subsection{Full model \label{sec::FullModel}}

Using the above descriptions of the rates of the three viral transmission pathways, we can formulate the following ODEs describing the evolution of the viral loads on fomite surfaces, and the hands and mucous membranes of individuals present/visiting the setting:

\begin{subequations}
\begin{align}
\begin{split}
\label{eqn::ODESurf}
\dv{\VirSurf[i](t)}{t}& = \sum_{k=1}^{\numhands}\Indic[k](t)\left(\transferrate[H_k]{S_i}-\transferrate[S_i]{H_k}\right) -\inact[S_i]-\resusprate[i]-\TraCleanSurf[i] \\
&+ \SDrate[S_i] + \sum_{k=1}^{\numhands}\Indic[k](t)\LDrate[k]{S_i}
\end{split} \\
\begin{split}
\label{eqn::ODEHand}
\dv{\VirHand[j](t)}{t}& =\Indic[j](t)\Bigg(\sum_{i=1}^{\numsurf}\left(\transferrate[S_i]{H_j}-\transferrate[H_j]{S_i}\right) + \left(\transferrate[M_j]{H_j} - \transferrate[H_j]{M_j}\right) -\inact[H_j] -\TraCleanHand[j]{j} \\
&+ \SDrate[H_j] + \sum_{k=1}^{\numhands}\Indic[k](t)\LDrate[k]{H_j} \Bigg)
\end{split} \\
\begin{split}
\label{eqn::ODEMem}
\dv{\VirMems[j](t)}{t}& =\Indic[j](t)\left(1 - \infec[j]\right) \Bigg(\left(\transferrate[H_j]{M_j}-\transferrate[M_j]{H_j}\right)-\inact[M_j]-\airrate[j]  +\frac{\VirAir}{\volair}\resp[j] \\
&+ \sum_{\substack{k=1 \\ k\neq j}}^{\numhands}\Indic[k](t)\LDrate[k]{M_j}  \Bigg)
\end{split}
\end{align}
\end{subequations}

The right-hand side of \cref{eqn::ODESurf} accounts for the fomite pathway (transfer from and to the hands of individuals while they are present in the setting; see \cref{eqn::IndicHeaviside}), the reduction of the viral load due to natural inactivation, resuspension from the surface and cleaning as well as the deposition of aerosol (small droplets, SD) and large droplets (LD) emitted by any individual in close proximity.

In \cref{eqn::ODEHand,eqn::ODEMem}, the right-hand side is non-zero only when individual $j$ is present in the setting (i.e., then the indicator function \cref{eqn::IndicHeaviside} is non-zero) and the model is able to track their interactions with the environment. Owing to this formulation, the viral loads on both their hands and mucous membranes remain constant when the individual is outside the setting (either before they enter it or after they leave) since the model is unaware of the individuals' whereabouts and interactions when they are not in the enclosed setting being modelled.

When individual $j$ is in the setting, the viral load on their hands follows \cref{eqn::ODEHand} which involves transfer to/from fomite surfaces and the individuals' own mucous membranes, natural deactivation of the virus and its removal through hand washing as well as the deposition onto the hands of small and large droplets.

In \cref{eqn::ODEMem} we assume that the viral load on the mucous membranes of infected individuals (for whom $\infec[j] = 1$) is so high that no interactions with the surroundings can either increase or reduce it. Therefore, the right-hand side of equation \cref{eqn::ODEMem} is multiplied by $(1-\infec[j])$ to 'freeze' the viral load on the mucosa of infected individuals at its initial level to avoid any artificial decrease in their viral loads due to shedding (also, the timescales considered here are assumed to be shorter than those leading to significant changes in the condition of infected individuals). For susceptible individuals, their mucous membrane viral load can be affected by self-inoculation with their hands, inhalation of virus-laden aerosol and the deposition of large droplets from infected individuals while in close contact with them, as well as viral inactivation and possible physical removal through respiratory activities before the virus can reach the lower respiratory tract.

Variations in the viral load within the total volume of air, $\volair$, of the setting are described by the following ODE:
\begin{equation}
\label{eqn::ODEAir}
\begin{split}
\dv{\VirAir}{t}&=\sum_{j=1}^{\numhands}\Indic[j]\left( \left(1-\fractionLD[j]\right)\airrate[j] - \frac{\VirAir}{\volair}\resp[j] \right) +\sum_{i=1}^{\numsurf}\left(\resusprate[i]-\SDrate[S_i]\right) \\
&-\sum_{i=1}^{\numhands}\SDrate[H_i]-\inact[air] - \vent
\end{split}
\end{equation}
where the only positive contributions are those due to aerosol produced by infected individuals while they are in the setting and possible resuspension of the virus from fomite surfaces. Aerosolised viral particles can be removed from the air by people inhaling them, via the physical removal of the virus due to ventilation with the rate $\vent = \ventrate\VirAir/\volair$ ($\ventrate$ is the ventilation air flow in \si{\metre\cubed \per \hour}), or through its natural inactivation. Note that, although the terms describing the deposition of small droplets, $\SDrate[i]$, and resuspension of the virus back into the air, $\resusprate[i]$, are included in this formulation, their magnitudes can be assumed to be negligible in comparison with the rate of natural inactivation in air and those of other transmission events. Therefore, their values are set to zero in all the simulations reported below.

We assume that start of the observation period the setting is uncontaminated, so that the initial conditions for the viral loads on fomite surfaces and in the air are zero:
\begin{subequations}
\begin{align}
& \VirSurf[i](t=0)=0, \qquad i=1,...,\numsurf \label{eqn::InitialSurf} \label{eqn::InitialCondSurf}\\
& \VirAir(t=0)=0
\label{eqn::InitialCondAir}
\end{align}
\end{subequations}
Similarly, susceptible individuals entering the setting are assumed to have had no prior exposure to the virus:
\begin{equation}
\label{eqn::InitialCondHandMemsSusceptible}
\VirHand[j](t=0)=\VirMems[j](t=0)=0, \qquad \forall j=1,...,\numhands, \; \infec[j]=0
\end{equation}
On the other hand, infected individuals carry significant viral loads $\VirMems[0]$ (see \Cref{tab:table_virus} of Appendix \ref{sec::AppendixA}) on their mucous membranes which are assumed to remain unchanged throughout the simulated period of $\tend$ hours:
\begin{equation}
\label{eqn::InitialCondMemsInfected}
\VirMems[j](t=0)=\VirMems[0], \qquad \forall j=1,...,\numhands, \; \infec[j] \neq 0
\end{equation}
Prior to entering the setting, infected individuals are also assumed not to have touched their faces and mucous membranes with typical frequencies but not any fomites, so that an equilibrium has been established between the viral loads on the individuals' mucous membranes and hands:
\begin{equation}
\label{eqn::InitialCondHandsInfected}
\VirHand[j](t=0) = \frac{\transferrate[M_j]{H_j} + \LDrate[j]{H_j}}{\transferconst[H_j]{M_j} + \inactconst[H_j]}, \qquad \forall j=1,...,\numhands, \; \infec[j] \neq 0
\end{equation}

\subsection{Generalised risk of infection}
Earlier, we defined the risk of infection in \cref{eqn::InfectionRisk} through the viral load accumulated on a susceptible individual’s mucous membranes as is commonly done in the literature \citep{Lei2018}. However, the ODE describing the viral load on mucous membranes in \cref{eqn::ODEMem} contains negative rates on the right-hand side describing the removal of viral particles through touch, respiratory activities and viral inactivation. To make sure the risk of infection is not underestimated, we will redefine it taking into account only the influx of viral particles to the mucous membranes while also considering separately the exposure through the aerosol, $\doseair[j]$, fomite, $\doseFom[j]$, and close contact (large droplet), $\doseLD[j]$, transmission pathways, respectively:
\begin{subequations}
\begin{align}
\label{eqn::DoseAir}
\doseair[j]&=\int_{\tstart}^{\tend} \Indic[j](t)\frac{\VirAir}{\volair} \resp dt
\\
\label{eqn::DoseFom}
\doseFom[j]&= \int_{\tstart}^{\tend} \Indic[j](t)\transferrate[H_j]{M_j} dt
\\
\label{eqn::DoseLD}
\doseLD[j]&= \int_{\tstart}^{\tend} \sum_{\substack{k=1 \\ k\neq j}}^{\numhands}\Indic[k](t)\LDrate[k]{M_j} dt
\end{align}
\end{subequations}

Now the risk of infection can be defined in the following way:
\begin{equation}
    \label{eqn::InfecRiskDose}
    \infecP[j]=1-\exp\left(-\frac{\left(\doseFom[j]+\doseair[j]+\doseLD[j]\right)}{\doserep[j]}\right)
\end{equation}
where we are explicitly assuming that the values of the dose response parameters $\doserep[j]$ associated with the different pathways are the same.

\subsection{Model implementation}
The mathematical model presented in this section was implemented in both MATLAB and Python. Refer to Appendix \ref{sec::AppendixB} for relevant implementation details.

\section{Model Demonstrations and Case Studies \label{sec::casestudies}}

To illustrate the generality and predictive ability of the model developed in \cref{sec::model}, we consider three scenarios of different complexity. Two of these, presented in \cref{sec::casestudy1}, are small-scale to enable a detailed analysis and interpretation of the effect of model parameters on the outputs, while the third scenario involves a larger group of people and explores how interactions within and between subgroups of a larger group affect individual risk of infection.

The values of model parameters describing the viral loads and shedding rate for infected individuals, viral transmission through droplets and between hands and fomites, survival of the virus on surfaces and in the air, etc., have been carefully collected from a number of literature sources and are summarised in \Cref{tab:table_virus} of Appendix \ref{sec::AppendixA}. Other parameter values, including reasonable assumptions about the values of some of the parameters described in the text below, are compiled in \Cref{tab:ltctable2} (for Case Studies 1 and 2) and \Cref{tab:ltctable1} for (for Case Study 3) in Appendix \ref{sec::AppendixA}. The models and related data have been made openly available on Github: https://github.com/Ishanki/VIRAS.

\subsection{Case Study 1: One-to-one Meeting \label{sec::casestudy1}}
In this case study we consider a small office or meeting room with air volume of \SI{40}{\m^3} (which corresponds to the floor area of \SI{16}{\m^2} with \SI{2.5}{\m} ceiling height), wherein two individuals one of whom is infected have a \SI{4}{\hour}-long face-to-face meeting (see \Cref{fig::SmallOfficeCaseStudies}a). While in the room, the individuals come into contact with three objects: a door handle, a desk and a document they are jointly working on. The individuals spend 90\% of the time they are together in \emph{close contact} so that the large respiratory droplets from the infected person can contaminate the desk, document and hands of the susceptible individual as well as directly deposit onto their mucous membranes. After the meeting, the infected individual leaves while the susceptible person remains in the room for another 4 hours. The room is ventilated at 1 air change per hour (ACH) or \SI{40}{\m^3/\hour} unless otherwise stated.

\SmallOfficeCaseStudies

In addition to ventilation, which continually dilutes the virus in the air, we also introduce the application of disinfecting agents on fomite surfaces and hands of both individuals as a way of viral inactivation. Note that, unlike hard surfaces, the shared document cannot be disinfected in this way. Hand and surface disinfection is performed either once after the meeting (i.e., at $t=\SI{4}{\hour}$), twice at $t=\SI{2}{\hour}$ and at $t=\SI{4}{\hour}$ or not at all. According to the mixed continuous-discrete formulation presented in \cref{sec::removal}, the viral loads on hands and fomites are diminished instantly upon disinfection by a factor $\cleanfracHand[j] = 10^{-\lrv[{\hand[j]}] }$ or $\cleanfracSurf[i] = 10^{-\lrv[{\surf[i]}]}$ according to \cref{eqn::handcleanrate} and \cref{eqn::fomitecleanrate}, respectively. In the simulations reported below these factors have the value of 100.

\Cref{fig::SmallOfficeLongContact} reports the simulated viral loads on the shared objects and in air, as well as the exposure of the susceptible individual to the virus. In the base case, when no surface and hand cleaning is performed, the viral loads on the shared objects increase throughout the duration of the meeting. The desk is characterised by the highest viral load, with the document accumulating approximately 2.5 times less virus, while the number of viral particles on the door handle at $t=\SI{4}{\hour}$ is 65 times lower than on the desk. 

Considering the traditional interpretation of the fomite transmission pathway, the difference in the viral loads should be directly proportional to the contact frequency, contact area and transferred fraction while being inversely proportional to an object's surface area according to \cref{eqn::RateXY}. Given the parameter values in \Cref{tab:ltctable2} and \Cref{tab:table_virus} of Appendix \ref{sec::AppendixA}, the rate of transmission from contaminated hands to the objects should be highest for the document, followed by the door handle and then by the desk. However, the viral loads in \Cref{fig::SmallOfficeLongContact}a contradict this conclusion. The difference between the expected and simulated fomite viral loads are explained by the contribution of large droplets which is the dominant route of contamination for the desk and the document that remain in close proximity of the infected individual for 90\% of the duration of the meeting.

When considering viral concentrations per unit of surface area, however, the most contaminated object is the document, followed by the door handle and then by the desk (data not shown), which is due to the differences in the fomite surface areas over which the viral particles are deposited.

After the meeting the viral loads on fomite surfaces begin to decrease as a result of natural viral deactivation in the absence of a source of contamination. The rate of decrease in each case is determined by the type of fomite surface, which affects the average survival time of the virus.

\Cref{fig::SmallOfficeLongContact}a also shows that disinfection of surfaces and hands at $t=\SI{4}{\hour}$ leads to a 100-fold decrease in the viral loads on hard surfaces (dashed lines). It is noteworthy that, even though the infected individual is no longer in the room after $t=\SI{4}{\hour}$, the viral loads begin to increase again. This is due to the continuing interactions of the susceptible individual with the still-contaminated document, which helps to spread the pathogen onto the disinfected surfaces. An additional cleaning event at $t=\SI{2}{\hour}$ results in a further decrease in the viral loads on fomite surfaces at later times, although the effect of an additional cleaning event is smaller.

As expected, air contamination remains independent of fomite disinfection events, as shown in \Cref{fig::SmallOfficeLongContact}b. The viral load increases throughout the meeting almost reaching a steady state when the rate of shedding of small contaminated droplets is approximately balanced by the combined rate of their removal by ventilation and natural viral inactivation in air. When the infected individual leaves the setting, the viral load in the air decreases exponentially owing to the latter two mechanisms.

It is of particular interest to observe the viral exposure of the susceptible person shown in \Cref{fig::SmallOfficeLongContact}c along with the contributions from the fomite, close contact, and aerosol pathways. The viral exposure continues to increase throughout the simulated $\SI{8}{\hour}$ duration, even after the infected person leaves the meeting, albeit at a slower rate. Indeed, for $t>\SI{4}{\hour}$ only one of the transmission pathways (close contact) is fully eliminated while the virus persists on fomite surfaces and in the air. It is also evident that the contribution of the fomite pathway is dominant throughout the $\SI{8}{\hour}$ period. As mentioned above, this is due primarily to shared surface contamination by large droplets, i.e., due to the fomite and large droplet pathways being intimately linked during prolonged periods of close proximity between infected and susceptible individuals. Furthermore, the aerosol pathway quickly becomes insignificant due to a combination of ventilation and a relatively short half-life of the virus in air. This is further illustrated in \Cref{fig::SmallOfficeLongContact}d, showing the diminishing with time of the relative contributions of the aerosol and close contact routes (i.e., direct transfer of large droplets to mucous membranes) compared to that of the fomite pathway.

\SmallOfficeLongContact

Owing to the importance of the fomite pathway in this case study the effect of surface and hand cleaning on viral exposure is expected to be significant, which is corroborated by the dashed and dash-dotted lines in \Cref{fig::SmallOfficeLongContact}c. Indeed, a single cleaning event at the end of the meeting leads to a 34\% reduction in the viral load accumulated by the susceptible person by the end of the working day, while two cleaning events give an overall 49\% reduction in viral exposure. The relative contributions of the different pathways also change: the contribution of fomite transmission decreases with enhanced cleaning and those due to close contact and aerosol increase (see \Cref{fig::SmallOfficeLongContact}d).

Such an important contribution of transmission via fomite surfaces warrants a further study into the effects of different mitigation strategies on the viral exposure of the susceptible individual. It is clear that increasing the frequency of hand and surface disinfection should lead to a significant reduction in exposure. However, it is also of interest to investigate how cleaning affects viral exposure in combination with other parameters such as (i) the close contact duration, (ii) face covering efficacy, and (iii) ventilation rate.

\Cref{fig::SmallOfficeHeatmap}a shows the synergistic effect on the final viral exposure of surface and hand disinfection (between 0 and 8 times during the $\SI{4}{\hour}$-long meeting) and variable close contact duration (i.e., the fraction of time when large droplets emitted by the infected person can land on mucous membranes of the susceptible one). Note that in this case the amount of large droplets landing on fomite surfaces and hands of the two people remains the same as in the results reported in \Cref{fig::SmallOfficeLongContact}. This can be interpreted as changing the relative positioning of the infected individuals which results in different amounts of droplets expelled directly into each other's face \cite{Zhang2020-dominantroutes}. The results show that the viral exposure decreases most rapidly with the increase in the number of cleaning events. When the number of cleaning events is greater than two (e.g., disinfection is performed more often than once every two hours) suppressing the direct transmission through large droplets can bring about a further reduction in the accumulated viral dose and hence the risk of infection.

\SmallOfficeHeatmaps

If large droplets can be captured at the source (e.g., by wearing a face covering) so that they are prevented from either landing on another person's mucous membranes or contaminating fomite surfaces, a very significant reduction in the viral exposure can be achieved even without hand and surface disinfection (\Cref{fig::SmallOfficeHeatmap}b) provided that the face covering captures the majority of large droplets. Together with hand and surface cleaning, a more significant reduction in viral exposure can be achieved by wearing efficient face coverings than in the case of reducing close contact only. This result showing how the large droplet and fomite pathways are intricately linked in these settings is an important observation from the model.

Lastly, \Cref{fig::SmallOfficeHeatmap}c shows the combined effect of cleaning and ventilation rate on viral exposure. The rate of ventilation has little effect on the viral exposure except at very low values (\textless 1\,ACH), while most of the reduction in this quantity in \Cref{fig::SmallOfficeHeatmap}c is due to hand and surface cleaning. This finding is not surprising when considering the low relative contribution of the aerosol pathway to the total viral exposure observed in \Cref{fig::SmallOfficeLongContact}d.

\subsection{Case Study 2: Working alongside in a small office \label{sec::casestudy2}}

The second case study uses the same environment as in Case Study 1 (\cref{sec::casestudy1}) but the two individuals are assumed to be working alongside each other (see \Cref{fig::SmallOfficeCaseStudies}b) for $\SI{4}{\hour}$, after which the infected person leaves the room and the susceptible one remains in the setting for a further $\SI{4}{\hour}$. Unlike in Case Study 1, the individuals spend only 5\% of the time they are together in the room in close proximity to each other and the other individual's desk, while each spending 50\% of the time handling the document. See \Cref{tab:ltctable2} of Appendix \ref{sec::AppendixA} for the other parameter values.

\SmallOfficeShortContact

The results of simulating the same scenarios as in Case Study 1 are presented in \Cref{fig::SmallOfficeShortContact}. The change in the way the individuals interact with each other and their environment leads to a different pattern of fomite contamination, shown in \Cref{fig::SmallOfficeShortContact}a, in which the document is less than half as contaminated as during prolonged close contact in Case Study 1 and the desk of the susceptible individual is over 40 times less contaminated than that of the infected person. There is no difference, however, in the contamination of the air (see \Cref{fig::SmallOfficeShortContact}b), since the aerosol pathway is not affected by the changes related to the fomite and close contact pathways.

Compared to Case Study 1, there is a nearly 7-fold decrease in the total viral exposure of the susceptible person at the end of the $\SI{8}{\hour}$-long working day (see \Cref{fig::SmallOfficeShortContact}c). This effect is explained by drastically reduced contributions of the fomite and close contact pathways due to reduced handling of shared objects (i.e., each other's desks and the document) as well as lower fomite contamination through the deposition of large droplets from the infected individual. This results in a different pattern of relative contributions of the three transmission pathways, as shown in \Cref{fig::SmallOfficeShortContact}d, in which the fomite pathway is still dominant overall, while the aerosol pathway dominates over relatively short times between approximately $\SI{0.25}{\hour}$ to $\SI{2}{\hour}$. However, it is important to emphasise here that the absolute viral exposure due to the aerosol pathway in this case is exactly the same as in Case Study 1 while the total exposure, and hence the risk of infection, are significantly lower than those predicted for Case Study 1.

Hand and fomite surface disinfection at $t = \SI{4}{\hour}$ leads to a significant reduction in the viral loads on surfaces as well as a 19\% decrease in the total exposure of the susceptible individual (\Cref{fig::SmallOfficeShortContact}a,c,d). The relative reduction in exposure is lower, however, than in Case Study 1 because of the lower contribution of the fomite transmission pathway into the total exposure. An additional cleaning event at $\SI{2}{\hour}$ leads to an overall decrease in viral exposure by 30\%.

The effect of combining mitigation strategies on the total viral exposure is reported in \Cref{fig::SmallOfficeHeatmap}d-f. First we note that further decreasing the fraction of time spent in close proximity from an already low value of 5\% (without changing the amount of large droplets contaminating fomites) leads to marginal improvements as compared to the enhanced hand and surface cleaning (\Cref{fig::SmallOfficeHeatmap}d). Capturing large droplets with a face covering of increasing efficacy has a more tangible effect on the viral exposure since this also reduces surface contamination (\Cref{fig::SmallOfficeHeatmap}e). The improvement compared to the results in \Cref{fig::SmallOfficeHeatmap}d is marginal, however, because of the relatively short amount of time the individuals spend in close proximity. \Cref{fig::SmallOfficeHeatmap}f shows that ventilation has a more pronounced influence on the total exposure than in Case Study 1 due to a higher relative importance of the aerosol pathway (although the absolute exposure is significantly lower, as mentioned above). At very low ventilation rates (\textless 1\,ACH), hand and surface cleaning become ineffective at reducing exposure since the aerosol pathway under such conditions becomes dominant. As the ventilation rate increases, the aerosol pathway becomes less important and disinfection becomes the main mitigation.

The two case studies, 1 and 2, illustrate the importance of tracking all known transmission pathways in enclosed settings in order to reveal which of them are dominant and over what time scales. It is also demonstrated how the model presented herein can be used to explore the effect of mitigation strategies on the viral exposure of susceptible individuals and thus pave the way to their optimal deployment in real life applications.

\subsection{Case Study 3: Graduate Student Office \label{sec::largetestcase}}
 To demonstrate the model's capability to track all transmission pathways while considering a large number of individuals and objects, a test case considering a graduate student office with one infected student and 38 susceptible students is used. It is based on the work and data of \citet{ZhangLi2018transmission}. While some of the features described in the original test case are retained, such as the categorisation of objects into private and public surfaces depending on personal or public use, several modifications have been made to demonstrate the model's wider applicability:
\begin{itemize}
\itemsep0em 
  \item Only desks and chairs are considered as private surfaces; all other personal belongings such as mugs, bags, etc. are ignored.
  \item Public surfaces include cabinet handles of 3 cabinets, printer, water dispenser, and door handle.
  \item Each cabinet is used by 13 students, where cabinet 1 is used by students 1 – 13, cabinet 2 by students 14 – 26, and cabinet 3 by students 27 – 39.
  \item For each public object, the surface area of the surface or part of the object with the highest touch frequency is considered. For example, the considered surfaces are: printer touch screen, water dispenser button, and cabinet handle.
  \item Average values are used for all contact frequencies, as opposed to individual contact frequencies based on individual behaviour, however, the model is also able to consider discrete touch events.
  \item Groups of friends, each containing 3 individuals, are considered, whereby the students in the group spend time in close contact with one another. It is assumed that students within each group do not come into close contact with each other’s private belongings.
\end{itemize}
Students in the office have been organised into 3 sets based on their interactions with the infected individual. Set 1 includes individuals who are in the friend group of and share a cabinet with the infected individual (Students 2 and 3). Set 2 includes individuals who share a cabinet with the infected individual, but are not in the same friend group as the infected individual (Students 4 – 13). Set 3 includes individuals who are neither in the friend group of the infected individual nor share a cabinet with them (Students 14 – 39). Furthermore, to highlight the impacts of different pathways, it is assumed that the infected individual spends 4 hours in the office while all susceptible students spend 8 hours. The parameters and respective values used in this test case are given in \Cref{tab:ltctable1} and \ref{tab:table_virus} of Appendix \ref{sec::AppendixA}. While data specific to SARS-CoV-2 has been used whenever available, note that some parameter values, such as transfer rates, have been obtained from data for the Influenza A virus \cite{ZhangLi2018transmission}, as these are currently unknown for SARS-CoV-2. 

\Cref{fig::Largecasegraphs} presents the results of this test case. It is evident from \Cref{fig::Largecasegraphs}a, which shows the relative exposures from each pathway, that the close contact pathway is the most dominant transmission pathway for students in Set 1, followed by the aerosol pathway. While the fomite pathway has the smallest contribution of approximately 7\%, it cannot be neglected in this instance. These results are further emphasised in \Cref{fig::Largecasegraphs}b, which shows the number of viral particles Student 2 (representative of Set 1) is exposed to via the different pathways over a 24-hour period. Noting that both Set 1 and Set 2 share a cabinet and therefore touch 'Cabinet handle 1', not shared by Set 3, exposure via the fomite pathway is negligible for individuals in both Sets 2 and 3. This indicates that shared public surfaces do not play a significant role in the fomite pathway in this instance. Further investigations reveal that large droplets deposited directly by the infected individual onto susceptible individuals' hands in Set 1 acts as the main contributor to the fomite pathway. As individuals in Set 2 do not come into close contact with the infected individual, this interaction between the close contact and fomite pathways is eliminated. 

\Cref{fig::Largecasegraphs}c shows the viral concentrations on public surfaces over a 24-hour period. Cabinet handle 1, shared by the infected individual and susceptible individuals in Sets 1 and 2, has the highest viral concentration throughout the 24-hour period. This is followed by the water dispenser button and printer touch screen. Although the water dispenser button has the highest touch frequency out of all public surfaces, cabinet handle 1 has the smallest surface area, potentially resulting in the higher viral concentration seen in \Cref{fig::Largecasegraphs}c. Cabinet handles 2 and 3 are not touched by the infected individual and therefore have negligible viral concentrations. Note that the viral concentrations in all public surfaces start to decrease rapidly after 4 hours as the infected individual exits at this time. Following the exit of all susceptible individuals at 8 hours, the viral concentrations continue to decrease via natural inactivation. \Cref{fig::Largecasegraphs}d. which shows the viral concentrations on Desk 2 and Chair 2 belonging to Student 2, representing Set 1, follows a similar trend over the 24-hour period. Despite the desk having a higher touch frequency when compared to the chair, the latter retains a higher viral concentration until approximately 18 hours, potentially due to its smaller surface area compared to the desk. 

\Largecasegraphs

Overall, this case study demonstrates the flexibility and adaptability of the model, which accommodates both continuous and discrete events and includes all three transmission pathways. The results shed light on how different pathways can dominate based on the extent of interaction between the susceptible and infected interactions, highlighting the importance of including all pathways, capturing human behaviour, and considering the interactions between pathways. Furthermore, the impacts of discrete events, such as the infected individual leaving the enclosed space, can be observed using this model as well.

\section{Conclusions}

We developed the model presented herein to quantitatively describe all the known transmission pathways of respiratory pathogens, such as SARS-CoV-2. The model formulation comprises, effectively, conservation equations tracking the virus with the level of detail sufficient to account for the expected (average) rates of transmission in a given setting: setting size and ventilation regime, objects present therein and their properties, the times of entry and exit of individuals and their interactions with other people and fomites. Like any model, the one developed here is based on a number of assumptions. First, we assume that the rate of dispersion of respiratory aerosols is high, so that any emitted small drops of mucus are instantaneously mixed with the air in the setting. This assumption is justified in most small to medium settings such as offices and classrooms. Second, we make assumptions about the fractions of large droplets that can be transmitted ballistically from one individual to mucous membranes of another or onto fomite surfaces when those are in close proximity to the source. These parameters are difficult to pinpoint in most cases, due to the distances between the source of large droplets and objects, as well as the orientation of the head of the infected individual constantly changing in most situations. Third, fomite transmission is described using either a continuous formulation using touch frequencies or through as a series of discrete events (the exact formulation can be chosen by the user depending on the context and available information). 

While more detailed representations of particular aspects such as exact airflow patterns or event-based simulations may provide more accurate predictions in some cases (particularly when scrutinising a transmission event \emph{post factum}), they can be too constrained to yield statistically relevant conclusions. Therefore, we believe that the level of approximation adopted in our model is sufficient for simulating the expected/average rates of transmission and risk of infection, as well as assessing the effect of mitigation measures, in general enclosed setting with less information available.

Our model offers new insights into the magnitude of exposure to the virus and the prevalent transmission routes under different conditions. The simulation results show that the nature of the enclosed setting and the intensity of person-to-person and person-fomite interactions therein play an important role in determining the relative contributions of the aerosol, close contact, and fomite transmission pathways. Thus, the results for Case Studies 1 and 2, focusing on two individuals sharing a small office, show that when individuals spend long periods of time in close proximity to each other and frequently touch shared objects the fomite pathway can be the dominant transmission route. The results also showed how transmission pathways can be related, particularly large droplet deposition and fomite pathways. In this case the contribution of the aerosol pathway may be significantly lower, and increasing the rate of ventilation (one of the most frequently recommended mitigation measures indoors) beyond a very modest 0.5\,ACH has a negligible effect on the risk of infection. A much more effective approach to reducing infection risk under such circumstances involves frequent cleaning of often-touched surfaces that are also in close proximity to more than one individual as well as wearing effective face covering to reduce the shedding of large droplets of mucus.

Our simulations also show that when close contact between individuals and between individuals and private objects of other people is minimised (as in Case Study 3), the overall exposure to the virus (and hence the risk of infection) are greatly reduced and transmission occurs predominantly through the aerosol pathway. This implies that the risk of infection under such conditions can be reduced through enhancing the rate of removal and/or inactivation of the pathogen in the air (e.g., through enhanced ventilation, air filtration or Ultraviolet germicidal irradiation).

These findings highlight the fact that the mechanisms of respiratory infection transmission can be drastically different depending on the nature of an indoor space, its occupancy level and how the occupants interact with each other and their environment. The time scales involved (e.g., the total length of stay in the setting and duration of close encounters with other people) are also important factors affecting the relative contributions of the different transmission pathways into total exposure. Based on our results, we can conclude that, while UK government COVID-19 advice undoubtedly lists valid mitigating interventions, their use should be tailored to particular enclosed environments based on a quantitative assessment of the different transmission pathways. This would allow not only a significant reduction in the infection risk in a particular setting to be achieved, but would also enable optimal deployment of mitigating interventions to minimise the associated monetary and non-monetary costs.

The model presented here provides a versatile means of simulating viral exposure at the level of a single enclosed space. However, it can also be used as part of a larger simulation involving multiple settings of different type visited by an individual throughout the day or over longer periods of time. Thus, this type of local, but high-fidelity, transmission model could play an important role as a building block in larger-scale epidemiological models.

\section{Acknowledgements}
The authors would like to express their gratitude to all those that participated in the Royal Society's Rapid Action in Modelling the Pandemic (RAMP), particularly those involved in Subgroup 4. The interactions and discussions in these meetings was always thought-provoking and inspired much of this work. We would also like to thank Dr Marco-Felipe King for his insights and encouragement during the initial model development. Finally, we are grateful to Dr Nan Zhang for providing us access to the data from their extensive work in \citet{NanZhangBehaviour}.

\newpage
\begin{appendices}
\section{Parameter details for SARS-CoV-2 Case Studies} \label{sec::AppendixA}

\subsection{Human behaviour and Objects}

The parameters used to describe human behaviour and interactions with objects are described in this section, for all three case studies. 

Table \ref{tab:ltctable2} contains the parameters specific to Case Studies 1 and 2, where two individuals interact with one another and fomites during a one-to-one meeting or while working alongside each other.

\newcommand{\ph}{\si{\per\hour}}
\newcommand{\cms}{\si{\centi\meter\squared}}
\begin{longtable}{ c c c c }
\caption{ Parameters used in case studies 1 and 2. The following subscripts have been used: doc - document, dsk - desk, dh - door handle, inf - infected individual, muc - mucosa, resp - respiration, susc - susceptible individual, vent - ventilation.}
\label{tab:ltctable2}\\
\hline \multicolumn{1}{c}{\textbf{Parameter}} & \multicolumn{1}{c}{\textbf{Value}} & \multicolumn{1}{c}{\textbf{Unit}} & \multicolumn{1}{c}{\textbf{Notes and Reference}}\\
\endfirsthead
%\hline Parameter & Value & Unit & Notes and References \\[8pt]
\hline 
\hline
$\Area[\doc]$ & $623.7$& \cms & Assumed \\
$\Area[\desk]$ & $6000$ & \cms & Assumed \\
$\Area[\doorh]$ & $65$& \cms & Assumed \\
$\Area[\muco]$ & $391.7$& \cms & Value calculated using data \cite{gizurarson2012anatomical,naumova2013oral,watsky1988comparison}\\
$\Area[\han]$ & $147.02$ & \cms & Area of both palms \cite{goker2017determination}\\
$\AreaCon[\han,]{\doc}$ & $36.8$ & \cms & Assumed \\
$\AreaCon[\han,]{\desk}$ & $73.5$ & \cms & Assumed \\
$\AreaCon[\han,]{\doorh}$ & $36.0$ & \cms & Assumed \\
$\AreaCon[\han,]{\muco}$ & $7.67$ & \cms & Value assumed for two fingertips \cite{sahmel2015estimation} \\
$\frequency[\han,]{\doc}$ & $5.0$ & \si{\per \hour} & Case Study 1, assumed \\
$\frequency[\han,]{\doc}$ & $2.5$ & \si{\per \hour} & Case Study 2, assumed \\
$\frequency[\han,]{\desk}$ & $20.0$ & \si{\per \hour} & Average value from \cite{ZhangLi2018transmission}\\
$\frequency[\han,]{\doorh}$ & $1$ & \si{\per \hour} & Assumed\\ 
$\frequency[\han,]{\muco}$ & $16$ & \si{\per \hour} & Average value from \cite{ZhangLi2018transmission}\\

$\resp[\susc \text{ or } \infectxt]$ & $0.39$ & \si{\metre\cubed \per \hour} & Average value\cite{airway2020levitan} \\
$\fractionCloseContact[\susc,]{\infectxt}$ & $0.9$ & - & Case Study 1, assumed \\
$\fractionCloseContact[\susc,]{\infectxt}$ & $0.05$ & - & Case Study 2, assumed \\
$\fractionCloseContact[\susc \text{ or } \infectxt,]{\doc}$ & $0.9$ & - & Case Study 1, assumed \\
$\fractionCloseContact[\susc \text{ or } \infectxt,]{\doc}$ & $0.5$ & - & Case Study 2, assumed \\
$\fractionCloseContact[\susc \text{ or } \infectxt,]{\desk}$ & $0.9$ & - & Case Study 1, assumed \\
$\fractionCloseContact[\susc,]{\desk1}$ & $0.9$ & - & Case Study 2, assumed \\
$\fractionCloseContact[\susc,]{\desk2}$ & $0.05$ & - & Case Study 2, assumed \\
$\fractionCloseContact[\infectxt,]{\desk1}$ & $0.05$ & - & Case Study 2, assumed \\
$\fractionCloseContact[\infectxt,]{\desk2}$ & $0.9$ & - & Case Study 2, assumed \\
$\fractionCloseContact[\susc \text{ or } \infectxt,]{\doorh}$ & $0.05$ & - & Assumed \\

$\fractionLDsurf[\doc]{}$ & $0.05$ & - & Assumed \\
$\fractionLDsurf[\desk]{}$ & $0.4$ & - & Assumed \\
$\fractionLDsurf[\doorh]{}$ & $0.006$ & - & Assumed \\
$\volair$ & $40$ & \si{\metre\cubed} & Assumed \\
$\ventrate$ & $40$ & \si{\metre\cubed \per \hour} & Unless otherwise stated \\

\hline

\end{longtable}

Table \ref{tab:ltctable1} contain those specific to Case Study 3, where infected and susceptible individuals in a large graduate student office are modelled.
\begin{longtable}{ c c c c }
\caption{ Parameters used in the large test case, based on the original test case by Zhang and Li \cite{ZhangLi2018transmission}. The following subscripts have been used: ch - cabinet handle, chr - chair, dsk - desk, dh - door handle, ps - printer touch screen, wdb - water dispenser button, inf - infected individual, group - specific to interactions within a group, muc - mucosa, resp - respiration, susc - susceptible individual, vent - ventilation.}
\label{tab:ltctable1}\\
\hline \multicolumn{1}{c}{\textbf{Parameter}} & \multicolumn{1}{c}{\textbf{Value}} & \multicolumn{1}{c}{\textbf{Unit}} & \multicolumn{1}{c}{\textbf{Notes and Reference}}\\
\endfirsthead
%\hline Parameter & Value & Unit & Notes and References \\[8pt]
\hline 
\hline
$A_{\desk}$ & 6000 & \cms & Desktop\cite{ZhangLi2018transmission} \\
$A_{\chair}$ & 4260 & \cms & Total chair area\cite{ZhangLi2018transmission} \\
$A_{\cabh}$ & 10 & \cms &  Handle area only\cite{ZhangLi2018transmission} \\
$A_{\printer}$ & 35 & \cms &  Touchscreen area only\cite{ZhangLi2018transmission} \\
$A_{\waterdisp}$ & 12 & \cms &  Button area only\cite{ZhangLi2018transmission} \\
$A_{\doorh}$ & 100 & \cms &  Handle area only\cite{ZhangLi2018transmission} \\

$A^{c}_{\desk}$ & 73.5 & \cms & Assumed \\
$A^{c}_{\chair}$ & 73.5 & \cms & Assumed \\
$A^{c}_{\cabh}$ & 7 & \cms & Assumed \\
$A^{c}_{\printer}$ & 17.5 & \cms & Assumed \\
$A^{c}_{\waterdisp}$ & 6 & \cms & Assumed \\
$A^{c}_{\doorh}$ & 70 & \cms & Assumed \\
$A^{c}_{\muco}$ & $7.67$ & \cms & Value assumed for two fingertips\cite{sahmel2015estimation} \\

$f_{\han,\desk}$ & $20.2$& \ph & Average value\cite{ZhangLi2018transmission} \\
$f_{\han,\chair}$& $8.9$ & \ph & Average value\cite{ZhangLi2018transmission} \\
$f_{\han,\cabh}$ & $0.14$& \ph & Average value\cite{ZhangLi2018transmission} \\
$f_{\han,\printer}$ & $0.28$& \ph & Average value\cite{ZhangLi2018transmission} \\
$f_{\han,\waterdisp}$ & $0.31$& \ph & Average value\cite{ZhangLi2018transmission} \\
$f_{\han,\doorh}$ & $0.05$& \ph & Average value\cite{ZhangLi2018transmission} \\
$f_{\han,\muco}$ & $16$& \ph & Average value\cite{ZhangLi2018transmission} \\

$\volair$ & $400$ & \si{\metre\cubed} & $ $ \cite{ZhangLi2018transmission} \\

$\ventrate$ & 400$ $ &\si{\metre\cubed \per \hour} & One air change per hour\cite{ZhangLi2018transmission} \\

$\fractionCloseContact[\susc \text{ or } \infectxt,]{\desk}$ & $0.8$ & - & Assumed for a 8-hr scenario \\
$\fractionCloseContact[\susc \text{ or } \infectxt,]{\chair}$ & $0.4$ & - & Assumed for a 8-hr scenario \\
$\fractionCloseContact[\susc \text{ or } \infectxt,]{\cabh}$ & $0.03$ & - & Assumed for a 8-hr scenario \\
$\fractionCloseContact[\susc \text{ or } \infectxt,]{\printer}$ & $0.01$ & - & Assumed for a 8-hr scenario \\
$\fractionCloseContact[\susc \text{ or } \infectxt,]{\waterdisp}$ & $0.005$ & - & Assumed for a 8-hr scenario \\
$\fractionCloseContact[\susc \text{ or } \infectxt,]{\doorh}$ & $0.005$ & - & Assumed for a 8-hr scenario \\
$\fractionCloseContact{group}$ & $0.125$ & - & Assumed for a 8-hr scenario \\

$\fractionLDsurf{\desk}$ & 0.07 &-& Assumed for a 8-hr scenario \\
$\fractionLDsurf{\chair}$ & 0.05 &-& Assumed for a 8-hr scenario \\
$\fractionLDsurf{\cabh}$ & 0.01 &-& Assumed for a 8-hr scenario \\
$\fractionLDsurf{\printer}$ & 0.05 &-& Assumed for a 8-hr scenario \\
$\fractionLDsurf{\waterdisp}$ & 0.05 &-& Assumed for a 8-hr scenario \\
$\fractionLDsurf{\doorh}$ & 0.01 &-& Assumed for a 8-hr scenario \\
$\fractionLDsurf{group,hand}$ & 0.01 &-& Assumed for a 8-hr scenario \\
$\fractionLDsurf{group,muc}$ & 0.005 &-& Assumed for a 8-hr scenario \\

\hline

\end{longtable}

\subsection{Physical properties}

Table \ref{tab:table_virus} contains parameters used in all three case studies, primarily those associated with viral transfer, transmission, and inactivation.
\begin{longtable}{ c c c c }
\caption{ Parameters associated with viral transmission and other related physical properties. The following subscripts have been used: inf - infected individual, muc - mucosa, resp - respiration, susc - susceptible individual, por - porous, npor - nonporous, ss - stainless steel
.}
\label{tab:table_virus}\\
\hline \multicolumn{1}{c}{\textbf{Parameter}} & \multicolumn{1}{c}{\textbf{Value}} & \multicolumn{1}{c}{\textbf{Unit}} & \multicolumn{1}{c}{\textbf{Notes and Reference}}\\
\endfirsthead
%\hline Parameter & Value & Unit & Notes and References \\[8pt]
\hline 
\hline

$\resusprate[j]$ & $0.39$ & \si{\metre\cubed \per \hour} & Average value for respiration\cite{airway2020levitan} \\
$\doserep[\susc]$ & $3.95 \times 10^{5}$ & $\text{\footnotesize{viral particles}}^{-1}$ & As deduced by \cite{zhang2020deducing} \\
$\lrv[i]$ & $2$ & - & Assumed based on other viruses \cite{Tuladhar2012}\\
$\lrv[\han]$ & $1.1$ & - & Assumed using data for Norwalk virus \cite{Liu2010}\\
$\lrv[\muco]$ & $0$ &-& Assumed \\
$\airrate[\infectxt]$ & $1.13015\times10^{7}$ & $\text{\footnotesize{viral particles}}^{-1}$ & Taken from \cite{Ma2020shedding} \\

$\TranRate[\han]{por}$ & $0.03$ & - & Value assumed from influenza data \cite{ZhangLi2018transmission} \\
$\TranRate[\han]{npor}$ & $0.07$ & -& Value assumed from influenza data \cite{ZhangLi2018transmission} \\ 
$\TranRate[\han]{ss}$ & $0.08$ & -& Value assumed from influenza data \cite{ZhangLi2018transmission} \\
$\TranRate[por]{\han}$ & $0.8$ &-&  Value assumed from influenza data \cite{ZhangLi2018transmission}\\
$\TranRate[npor]{\han}$ & $0.12$ & -& Value assumed from influenza data \cite{ZhangLi2018transmission} \\
$\TranRate[ss]{\han}$ & $0.16$ & -&  Value assumed from influenza data \cite{ZhangLi2018transmission} \\
$\TranRate[\muco]{\han}$ & $0.5$ &-& Value assumed from influenza data \cite{ZhangLi2018transmission} \\

$\VirMems[0]$ & $4.0\times10^{6}$ & viral particles & Assumed\\

$\cleaneffecSurf$ & $1$ &-& Assumed for all surfaces\\
$\fractionLD$ & $0.5$ &-& Assumed \\

$\halflife[\air]$ & $1.1$&\si{\hour} & Value obtained from \cite{van2020aerosol}\\
$\halflife[por]$ & $3.46$ & \si{\hour} & Value obtained from \cite{van2020aerosol}\\
$\halflife[npor]$ & $6.81$  & \si{\hour} & Value obtained from \cite{van2020aerosol}\\
$\halflife[ss]$ & $5.63 $  & \si{\hour} & Value obtained from \cite{van2020aerosol}\\
$\halflife[\han]$ & $3.5$  & \si{\hour} & Sars-CoV-2 half-life on swine skin at \SI{22}{\degreeCelsius} \cite{XiaoRoleofFomites} \\

\hline

\end{longtable}

\section{Implementation notes} \label{sec::AppendixB}
Both the continuous and discrete model formulations of the model have been implemented in MATLAB and Python. The resulting system of ODEs is solved numerically using \texttt{ode15s} solver in MATLAB or SciPy's \texttt{odeint()} in the Python implementation.

To avoid discontinuities on the right-hand sides of \cref{eqn::ODESurf,eqn::ODEHand,eqn::ODEMem,eqn::ODEAir} due to the appearance of both the Dirak's delta-function in \cref{eqn::handcleanrate} and \cref{eqn::fomitecleanrate} and indicator functions \cref{eqn::IndicHeaviside} we replace them with continuous approximations using a triangular function to approximate Dirak's delta function:
\begin{equation}
    \tridelta(t)=\left\lbrace \begin{array}{cc}
      0,   & |t|>\trieps \\
\frac{1}{\trieps}\left(1-\frac{|t|}{\trieps}\right), & |t|\leq \trieps
    \end{array}\right.
\end{equation}
and the following piecewise linear approximation of the Heaviside function:
\begin{equation}
    \label{eqn::approxHeaviside}
    \triHeaviside(t)=\left\lbrace \begin{array}{cc}
        0,   & t<0 \\
        x/\trieps, & 0\leq t \leq \trieps \\
        1,   & t>\trieps \\
    \end{array}\right.
\end{equation}

The discontinuous indicator function \cref{eqn::IndicHeaviside} is then replaced in the code by the following trapezoidal function:
\begin{equation}
    \label{eqn::ApproxIndic}
    \ApproxIndic[j](t)=\triHeaviside(t-\tin[j]) \, \triHeaviside(\tin[j]+\tdur[j]+\trieps-t), \ \ j = 1,\ldots,\numhands
\end{equation}

One, however, must ensure that the parameter $\trieps$ is much smaller than the shortest duration of occupancy:
\begin{equation}
    \trieps << \max_{j = 1,\ldots,\numhands} \tdur[j]
\end{equation}
and the maximum integration step size used by \texttt{ode15s} or \texttt{odeint()} is significantly less than $\trieps$.

\end{appendices}

\bibliography{References}

\end{document}